\documentclass[preprint,12pt]{elsarticle}
\usepackage{graphicx} 
\graphicspath{ {./figures/} }
\usepackage[rightcaption]{sidecap}
\usepackage{caption}
\usepackage{subcaption}
\usepackage{wrapfig}
\usepackage{float}
\usepackage{imakeidx}
\usepackage{comment}
\usepackage{commath}
\usepackage[titletoc]{appendix}
\usepackage{graphicx}
\usepackage{adjustbox}   
\usepackage{resizegather}
\usepackage[english]{babel}
\usepackage{tabu}
\usepackage{booktabs}
\usepackage{xfrac}
\usepackage{tabularx}
\usepackage{amssymb}
\usepackage{amsmath}
\usepackage{commath}
\usepackage{graphicx,bm}
\usepackage{verbatim}
\usepackage{lscape}
\usepackage{relsize}
\usepackage{enumitem}
\usepackage{textcomp}
\usepackage{breqn}
\usepackage{makecell}
\usepackage{longtable}
\usepackage{multirow}
\usepackage{doi}
\usepackage{fancyhdr}
\usepackage{algorithm}
\usepackage{algpseudocode}
\usepackage{setspace}
\usepackage{footnote}
\PassOptionsToPackage{hyphens}{url}
\usepackage{hyperref}
\hypersetup{colorlinks,linkcolor={blue},citecolor={blue},urlcolor={blue}}
\usepackage[numbers]{natbib}
\usepackage{mathtools}
\usepackage[disable]{todonotes}
\usepackage[framemethod=tikz]{mdframed}
\usepackage{booktabs,xcolor,siunitx}
\usepackage{soul}
\usepackage{cleveref}
\usepackage{tikz}
\usepackage{rotating}
\usetikzlibrary{shapes,snakes,automata,positioning}
\usepackage{tabularx}
\newcolumntype{Y}{>{\centering\arraybackslash}X}
\newcounter{mysubequations}

\newcommand{\f}[2]{\frac{#1}{#2}}
\newcommand{\mb}[1]{\mathbf{#1}}

\DeclareMathAlphabet\mathbfcal{OMS}{cmsy}{b}{n}
\renewcommand{\d}{\mathop{}\!\mathrm{d}} 
\newcommand{\p}{\partial}
\newcommand{\vecop}{\operatorname{vec}}

\definecolor{rev1}{HTML}{FF999A} 
\definecolor{rev2}{HTML}{F3F298} 
\definecolor{rev3}{HTML}{B2E0AE} 
\definecolor{other}{HTML}{C8C7FF} 

\title{DKFNet: Differentiable Kalman Filter for \\
Field Inversion and Machine Learning} 

\author{Yuan Wu\footnote{Graduate Student, University of Tennessee, Knoxville and Rice University}}
\author{Sicheng He\footnote{Assistant Professor,  sicheng@utk.edu, University of Tennessee, Knoxville}}

\date{\today}

\begin{document}

\begin{abstract}
The Kalman filter is a fundamental tool for state estimation in dynamical systems.
While originally developed for linear Gaussian settings, it has been extended to nonlinear problems through approaches such as the extended and unscented Kalman filters.
Despite its broad use, a persistent limitation is that the underlying approximate model is fixed, which can lead to significant deviations from the true system dynamics.
To address this limitation, we introduce the differentiable Kalman filter (DKF), an adjoint-based two-level optimization framework designed to reduce the mismatch between approximate and true dynamics.
Within this framework, a field inversion step first uncovers the discrepancy, after which a closure model is trained to capture the discovered dynamics, allowing the filter to adapt flexibly and scale efficiently.
We illustrate the capabilities of the DKF using two representative examples: a rocket dynamics model and the Allen–Cahn boundary value problem.
In both cases, and across a range of noise levels, the DKF consistently reduces state reconstruction error by at least 90\% compared to the classical Kalman filter, while also maintaining robust uncertainty quantification.
These results demonstrate that the DKF not only improves estimation accuracy by large margins but also enhances interpretability and scalability, offering a principled pathway for combining data assimilation with modern machine learning.
\end{abstract}

\maketitle{}

\section{Introduction}

The Kalman filter is a recursive algorithm used to estimate the hidden state of a linear dynamical system from noisy observations, assuming Gaussian noise and known system models. 
It operates by predicting the next state using the system dynamics and then updating that prediction using new measurements, balancing model uncertainty and measurement noise through a dynamically computed gain. 
Originally developed for aerospace navigation, the Kalman filter is now widely applied in fields such as autonomous vehicle localization ~\cite{Rezaei2007a,Liu2021a,Lu2022a,Chalvatzaras2023a}, satellite tracking ~\cite{Teixeira2008a,Zahaby2009a,Yang2016a,Pei2022a}, financial time series forecasting ~\cite{Huang2013a,Bao2014a,Khashei2020a}, and robotic control ~\cite{Rodriguez1987a,Gautier2001a,Du2015a,Martin2018a,Lee2020a}, where it enables accurate real-time estimation even with incomplete or noisy sensor data.

Multiple variants of the classic Kalman filter have been proposed in recent years to address the challenges encountered in its application.
The popular ones include the extended Kalman filter (EKF), unscented Kalman filter (UKF), ensemble Kalman filter (EnKF), and particle filter (PF).
The EKF linearizes nonlinear systems using a first-order Taylor expansion, but its accuracy drops when nonlinearities are strong \cite{Humpherys2012a,Kao2004a,Branicki2012a}. 
The UKF overcomes this by propagating sigma points through the nonlinear system, providing more accurate estimates in highly nonlinear cases, though at increased computational cost for large state spaces \cite{Gyrgy2014a,Kang2017a,Ijaz2008a}.
The EnKF uses an ensemble of forecasts for efficient state estimation in high-dimensional and nonlinear systems, but may need large ensemble sizes to accurately represent uncertainty \cite{Evensen2003a,Grooms2014a,Harlim2014a,Xie2025a,Sebacher2013a}. 
The PF, in contrast, approximates the full posterior distribution with weighted particles, making it highly flexible for nonlinear and non-Gaussian models, but often suffers from sample degeneracy and high computational demand, especially as the state dimension grows \cite{DelMoral2015a}.

In recent years, there has been a surge of interest in extending the Kalman filter with differentiable programming techniques, enabling end-to-end learning of both system and measurement models directly from data \cite{Kloss2021a,Liu2023a,Shen2025a}. 
Such approaches --- which we broadly refer to here as Kalman filters using differentiable methods --- include variants such as the differentiable extended Kalman filter (DEKF), the differentiable unscented Kalman filter (DUKF), differentiable ensemble Kalman filter (DEnKF). 
For a comparison of current works and our new results, see \Cref{tab:method_comparison}.
The DEKF proposed by \cite{Piga2021a} allows both the system dynamics and observation models to be identified directly from tactile sensor data, without requiring hand-crafted models.
Unlike classic EKF approaches that require hand-crafted physical models, differentiable EKF allows the entire state estimation pipeline to be trained via gradient-based optimization. 
This framework was demonstrated on the iCub humanoid robot, an open-source research platform developed for studying embodied AI, where it successfully tracked the position and velocity of sliding objects using only tactile observations, achieving a high degree of accuracy. 
The DEKF thus provides a flexible, data-driven approach for state estimation in scenarios where conventional modeling is difficult or inaccurate.
However, its effectiveness can be hindered by practical challenges such as sensitivity to model mismatch and limited robustness when applied to highly nonlinear or noisy systems.

In \cite{Revach2022a}, the DUKF enables end-to-end learning of both system and uncertainty models directly from data.
Unlike unstructured deep learning methods, DUKF retains the model structure of the classic UKF, which improves interpretability—that is, the ability to understand and trace how state estimates are computed based on explicit filtering steps and model assumptions.
By embedding learning within the established UKF framework, DUKF provides both adaptability and a transparent, physically meaningful estimation process.
However, its applicability remains limited, as DUKF does not naturally scale to high-dimensional problems and cannot support gradient-based end-to-end training, which constrains its use in modern large-scale learning scenarios.

DEnKF has further broadened the applicability of these techniques to high-dimensional problems, such as soft robots dynamics estimation \cite{Liu2023b} and hybrid Bayesian experimental design for model discrepancy calibration \cite{Chen2022a}.
It leverages an ensemble of particles to efficiently represent uncertainty and propagate information through high-dimensional state spaces. 
By enabling end-to-end training, it can adapt both the underlying dynamics and observation models while retaining the statistical interpretability and parallelizability of classic EnKF frameworks, making it especially suitable for complex, data-rich environments.
However, its reliance on linear updates can limit accuracy in strongly nonlinear systems, and the need for large ensembles increases computational cost, which constrains its scalability in practice.

Beyond the Kalman filter family, differentiable filter has also advanced particle-based methods. 
Differentiable particle filter (DPF) integrates neural networks with algorithmic priors for end-to-end optimization of the entire particle-filtering process \cite{Corenflos2021a}. 
The proposal mixture neural network for differentiable particle filters (PropMixNN) uses neural networks to learn proposal distributions in particle filters, reducing estimation errors in highly nonlinear systems \cite{Cox2024a}. 
In the context of physics-constrained filter,~\cite{Yin2024a,Jiang2022a} proposed a physics-constrained dynamic mode decomposition (PCDMD) framework that integrates Kalman filter with physical residuals, leading to improved dynamic mode decomposition based forecasts under noisy or imperfect models.
\begin{table}[htbp]
\centering
\caption{Comparison of classic differentiable filtering and our methods.}
\scriptsize                      
\setlength{\tabcolsep}{3pt}      
\renewcommand{\arraystretch}{1.2}
\begin{adjustbox}{width=\columnwidth,center}  
\begin{tabularx}{\columnwidth}{lYYYYY}
\hline
Method & \makecell{Learning/\\Adaptation}
& \makecell{End-to-end\\Gradient}
& \makecell{Nonlinearity\\Handling}
& \makecell{High-Dim\\Capable}
& \makecell{Key\\Publications} \\
\hline
DEKF  & $\times$          & $\times$ & $\times$       & $\times$       
& \cite{Piga2021a}                        \\
DUKF  & \checkmark          & $\times$ & $\times$       & $\times$       
& \cite{Revach2022a}                    \\
DEnKF & \checkmark          & $\times$ & \checkmark     & $\times$       
& \cite{Liu2023b,Chen2022a}                     \\
DPF   & \checkmark       & $\times$ & \checkmark     & $\triangle$    
& \cite{Corenflos2021a,Cox2024a} \\
\textbf{DKF}
& \textbf{\checkmark} & \textbf{\checkmark}
& \textbf{\checkmark} & \textbf{\checkmark}
& \textbf{This work}                    \\
\hline
\end{tabularx}
\end{adjustbox}
\label{tab:method_comparison}
\vspace{1mm}

\noindent
\textbf{Symbol meanings:}
$\checkmark$: Supported;\ 
$\times$: Not supported;\ 
$\triangle$: Partially supported.  
\end{table}

These works illustrate how differentiable filtering ideas extend beyond Kalman-type algorithms to particle-based and physics-constrained formulations. 
Although conceptually distinct from sequential filtering, the field inversion and machine learning (FIML) framework introduced by \cite{Parish2016a} is closely related in spirit: both aim to augment imperfect physics models with data-driven corrections within a gradient-based optimization framework. 
In FIML, this is achieved by directly inferring spatially distributed functional corrections for computational physics models, rather than simply tuning model parameters.
This approach systematically addresses model-form errors by first performing field inversion using observational or high-fidelity data to obtain corrective terms, and then applying machine learning techniques to reconstruct these corrections as functions of relevant variables. 
The resulting corrected models are then used to augment predictive simulations, significantly improving model accuracy and providing quantified model-form uncertainty.
However, classic FIML also has important limitations. 
It is typically performed as a one-time inversion on a fixed dataset, so the learned corrections remain unchanged during prediction and cannot be adapted as new observations become available. 
Moreover, its effectiveness depends strongly on access to high-fidelity data and accurate specification of observational covariances, and the inversion process itself can be computationally expensive in high-dimensional settings. 
Most importantly, FIML focuses on correcting model-form errors in a batch setting but does not integrate with sequential state estimation and uncertainty propagation, leaving a gap that motivates our DKF framework.

To overcome these limitations, our approach formulates the state transition operator as a learnable parameter, refined through field inversion by minimizing the mismatch between predicted and observed data using end-to-end gradient-based optimization. 
After optimization, a deep neural network is trained to capture the improved dynamics, allowing for flexible and accurate model correction. 
Compared with classic FIML, which performs one-time batch corrections detached from sequential filtering, our method integrates model correction directly into the DKF pipeline, enabling both states and dynamics to be jointly adapted as new data become available. 
In contrast to regular DKF approaches that rely solely on automatic differentiation, we incorporate analytically derived gradients, which improve efficiency, stability, and interpretability while preserving mathematical rigor. 
This systematic design bridges the gap between data assimilation and modern machine learning, providing a principled and scalable way to discover and adapt the underlying dynamics within the DKF pipeline.

The rest of the paper is organized as follows.  
We first review the classic Kalman filter framework in \Cref{sec:KF} to motivate the DKF. 
Next, in \Cref{sec:DKF}, we present the main algorithm, including the DKF and closure model using machine learning.  
Finally, we apply the algorithm to two data-driven models in \Cref{sec:experiments}.  

\section{Kalman filter }
\label{sec:KF}
The Kalman filter, named after Rudolf E. Kálmán, is an efficient recursive algorithm used to estimate the state of a linear dynamic system from a series of noisy measurements \cite{Kalman1960a}. 
The ground truth dynamics and measurement are defined by
\begin{equation}
\label{eq:dynamics and measurement}
\begin{aligned}
(\text{Dynamics}) \quad \mathbf{x}_{k} &= \mathbf{F}_k \mathbf{x}_{k-1} + \mathbf{B}_k \mathbf{u}_k + \mb{f}_k + \mathbf{w}_k, \\
(\text{Measurement}) \quad \mathbf{z}_{k} &= \mathbf{H}_k \mathbf{x}_{k} + \mathbf{v}_k,
\end{aligned}
\end{equation}
where $\mb{w}_k\in\mathbb{R}^{n_{{x}}}\sim \mathcal{N}(\mb{0}, \mb{Q}_k)$, and  $v_k\in\mathbb{R}^{n_{{z}}}\sim \mathcal{N}(\mb{0}, \mb{R}_k)$ are the external noise.
$\mb{F}_k\in\mathbb{R}^{n_{{x}}\times n_{{x}}}$, $\mb{H}_k\in\mathbb{R}^{n_{{z}}\times n_{{x}}}$, $\mb{B}_k\in\mathbb{R}^{n_{{x}}\times n_{{u}}}$,
$\mathbf{u}_k \in \mathbb{R}^{n_u}$, 
$\mb{f}_k \in \mathbb{R}^{n_x}$, 
$\mb{Q}_k\in\mathbb{R}^{n_{{x}}\times n_{{x}}}$ and $\mb{R}_k\in\mathbb{R}^{n_{{z}}\times n_{{z}}}$ denote state transition matrix, observation matrix, control input matrix, control vector, deterministic forcing, covariance matrix of the process noise and the observation noise.
The model above describe a linear case. 
The framework can also be extended to nonlinear dynamics and observation operators.

However, when numerical models are flawed and deviate from the true system, the state--space relations can only be represented approximately:
\begin{equation}
\label{eq:imperfect_model}
\begin{aligned}
(\text{Model dynamics}) \quad 
\hat{\mathbf{x}}_{k} &= \hat{\mathbf{F}}_k \hat{\mathbf{x}}_{k-1} 
                     + \hat{\mathbf{B}}_k \hat{\mathbf{u}}_k 
                     + \hat{\mathbf{f}}_k + \hat{\mathbf{w}}_k, \\
(\text{Observation}) \quad 
\hat{\mathbf{z}}_{k} &= \hat{\mathbf{H}}_k \hat{\mathbf{x}}_{k} 
                     + \hat{\mathbf{v}}_k,
\end{aligned}
\end{equation}
where the ``hatted'' operators denote the imperfect model counterparts of the true system matrices:  
$\hat{\mathbf{F}}_k \in \mathbb{R}^{n_x \times n_x}$ is the approximate state transition matrix,  
$\hat{\mathbf{B}}_k \in \mathbb{R}^{n_x \times n_u}$ the control input matrix,  
$\hat{\mathbf{H}}_k \in \mathbb{R}^{n_z \times n_x}$ the observation operator, and  
$\hat{\mathbf{f}}_k \in \mathbb{R}^{n_x}$ the deterministic forcing.  
The vectors $\hat{\mathbf{x}}_k \in \mathbb{R}^{n_x}$ and $\hat{\mathbf{z}}_k \in \mathbb{R}^{n_z}$ represent the model-predicted state and measurement, respectively.  
The process and measurement noises are given by  
$\hat{\mathbf{w}}_k \in \mathbb{R}^{n_x} \sim \mathcal{N}(\mathbf{0}, \hat{\mathbf{Q}}_k)$ and  
$\hat{\mathbf{v}}_k \in \mathbb{R}^{n_z} \sim \mathcal{N}(\mathbf{0}, \hat{\mathbf{R}}_k)$,  
with covariance matrices  
$\hat{\mathbf{Q}}_k \in \mathbb{R}^{n_x \times n_x}$ and  
$\hat{\mathbf{R}}_k \in \mathbb{R}^{n_z \times n_z}$.

At each step $k$, we define the residuals of the state and covariance as
\begin{equation}
\label{eq:KF_residual_step}
\mathbf{r}_{k} =
\begin{bmatrix}
\mathbf{r}_{\mathbf{x},k} \\
\mathbf{R}_{\mathbf{P},k}
\end{bmatrix}
=
\begin{bmatrix}
{\hat{\mathbf{x}}}_{k} - (\mathbf{I} - \mathbf{K}_k \hat{\mathbf{H}}_k)(\hat{\mathbf{F}}_k \hat{\mathbf{x}}_{k-1} + \hat{\mathbf{B}}_k \hat{\mathbf{u}}_k + \hat{\mathbf{f}}_k) - \mathbf{K}_k \hat{\mathbf{z}}_k \\
\mathbf{P}_{k} - (\mathbf{I} - \mathbf{K}_k \hat{\mathbf{H}}_k)(\hat{\mathbf{F}}_k \mathbf{P}_{k-1} \hat{\mathbf{F}}_k^\top + \hat{\mathbf{Q}}_k)
\end{bmatrix},
\end{equation}
stacking over all time steps yields the global residual vector
\begin{equation}
\label{eq:KF_residual}
\mathbf{r}(\mathbf{x}, \mathbf{d})
=
\begin{bmatrix}
\mathbf{r}_{1} \\
\mathbf{r}_{2} \\
\vdots \\
\mathbf{r}_{n_t}
\end{bmatrix}
= \mathbf{0},
\end{equation}
where $\mathbf{x} \coloneqq [\hat{\mathbf{x}}_1,\ldots,\hat{\mathbf{x}}_{n_t}]$ denotes the stacked state vector over the full time horizon, 
$\mathbf{d} \coloneqq [\hat{\mathbf{F}}_1,\ldots,\hat{\mathbf{F}}_{n_t}]$ collects the dynamics operators, 
and $\mathbf{K}_{k}$ is given by
\begin{equation}
\mathbf{K}_k = (\hat{\mathbf{F}}_k \mathbf{P}_{k-1} \hat{\mathbf{F}}_k^\intercal + \hat{\mathbf{Q}}_k)\, \hat{\mathbf{H}}_k^\intercal 
\left(\hat{\mathbf{H}}_k (\hat{\mathbf{F}}_k \mathbf{P}_{k-1} \hat{\mathbf{F}}_k^\intercal + \hat{\mathbf{Q}}_k) \hat{\mathbf{H}}_k^\intercal + \hat{\mathbf{R}}_k \right)^{-1}.
\end{equation}
details are provided in~\ref{sec:Kalman filter derivation}.
In the following~\Cref{sec:DKF}, we introduce DKF framework.

\section{DKF framework optimization}
\label{sec:DKF}

The Kalman filter is a useful tool of data assimilation, which combines the model and the observation to give the best prediction into the future.
However, the Kalman filter alters the dynamics of the system indirectly, not affecting the operator, ${\hat{\mb{F}}}_k$.
In practice, the physics models are usually inaccurate and many times miss critical physics happening in the real world, in other words, ${\hat{\mb{F}}}_k$ is merely an approximation of the real physics and can be further improved.
We propose to use the DKF to help improve the modeling of the underlying physics via gradient-based optimization in~\Cref{subsec:Differentiable Kalman filter}.
Then, a closure model is trained to discover the underlying dynamics using the optimization solution in~\Cref{subsec:Closure model}.

\subsection{Differentiable Kalman filter (DKF)}
\label{subsec:Differentiable Kalman filter}
The DKF formulates the entire Kalman filter process as a differentiable computation, enabling gradients to be computed with respect to model parameters. 
This allows for efficient field inversion of the dynamics, where the transition operator is optimized to minimize prediction errors, and for rigorous sensitivity computation using adjoint-based methods.
We present the following optimization problem:
\begin{equation}
\label{eq:opt}
\begin{aligned}
\min \,\, & ||\mb{e}||_2  \\
\text{w.r.t.}\,\, & \hat{\mb{F}}_1, \ldots, \hat{\mb{F}}_{n_t},
\end{aligned}
\end{equation}
where the deviation from the true observations serves as the measure of model accuracy, 
captured by the stacked residual vector
\begin{equation}
\label{eq:residual_vector}
\mathbf{e} = 
\begin{bmatrix}
\left(\hat{\mathbf{H}}_1 \hat{\mathbf{x}}_1 - \hat{\mathbf{z}}_1\right)^\intercal &
\ldots &
\left(\hat{\mathbf{H}}_{n_t} \hat{\mathbf{x}}_{n_t} - \hat{\mathbf{z}}_{n_t}\right)^\intercal 
\end{bmatrix}^\intercal,
\end{equation}
with $\hat{\mathbf{H}}_k$ the observation matrix, $\hat{\mathbf{x}}_k$ the estimated state, and $\hat{\mathbf{z}}_k$ the observation at time step $k$, for $n_t$ total time steps. 
The dynamics is implicitly enforced with the uncertainty treated as ``frozen''. 
Details are in the following~\Cref{Subsubsec:Field inv.} and~\Cref{subsubsec:Sensitivity computation}.

\subsubsection{Field inversion of the dynamics}
\label{Subsubsec:Field inv.}
We treat the filter's prediction step as a field inversion problem, holding the uncertainty fixed while we adjust the transition operator.  
Beginning with an initial guess~\(\hat{\mathbf{F}}_k^{(0)}\), we use gradient-based optimization to find a refined operator~\(\hat{\mathbf{F}}_k^{(*)}\) that minimizes the mismatch between our predicted measurements and the actual observations.  
In practice, this means running the standard Kalman predict–update cycles, but allowing the state‐transition matrix to adapt so that the filter's forecast residuals are as small as possible.

We begin by embedding the transition operator directly in the Kalman predict-update loop and then perform gradient‐based optimization on those matrices to minimize the forecast residuals.
By iterating this procedure, we discover a refined, state‐dependent approximation of the true dynamics.  
Finally, we replace the original, hand‐tuned operator with the optimized matrices, closing the loop and yielding more accurate predictions and more reliable uncertainty estimates in subsequent filtering steps.


\subsubsection{Sensitivity computation}
\label{subsubsec:Sensitivity computation}
The field inversion procedure described in the previous~\Cref{Subsubsec:Field inv.} relies on optimizing the state transition matrices $\hat{\mb{F}}_k$ so that the predicted observations from the Kalman filter best match the actual measurements. 
This optimization is achieved through a gradient-based approach, which in turn depends critically on the ability to efficiently compute the sensitivity of the loss function with respect to the design variables.

Specifically, we formulate the prediction error across all time steps as the loss function Eq.~\eqref{eq:opt}, and our design variable vector $\mb{d}$ is constructed by $[\hat{\mb{F}}_1, \hat{\mb{F}}_2, \ldots , \hat{\mb{F}}_{n_{t}}]$. 
To perform the optimization in Eq.~\eqref{eq:opt}, we require the gradient $\d f / \d \mb{d}$. 
Because the Kalman filter involves a sequence of recursive operations—each depending on previous estimates, covariances, and transition matrices—directly computing this gradient is computationally prohibitive, especially in high dimensions. Instead, we adopt the adjoint method, which is both efficient and analytically tractable for large-scale systems.

Following the methodology introduced by \cite{Jameson1988a},
We consider a system governed by a set of residual equations 
which we have defined in Eq.\eqref{eq:KF_residual}.
Our goal is to minimize a function of interest $f(\mathbf{x}, \mathbf{d})$. 
The adjoint method provides an efficient way to compute the sensitivity 
$\mathrm{d} f / \mathrm{d} \mathbf{d}$ via
\begin{equation}
\label{eq:adjoint}
\begin{aligned}
\left(\frac{\partial \mathbf{r}}{\partial \mathbf{x}}\right)^\intercal \boldsymbol{\psi} &= \frac{\partial f}{\partial \mathbf{x}}, \\
\frac{\mathrm{d} f}{\mathrm{d} \mathbf{d}} &= \frac{\partial f}{\partial \mathbf{d}} - \boldsymbol{\psi}^\intercal \frac{\partial \mathbf{r}}{\partial \mathbf{d}}.
\end{aligned}
\end{equation}

Incorporating the equation, we have
\begin{equation}
\renewcommand{\arraystretch}{1.2}
\begin{bmatrix}
\mb{A}_{1,1} & 0      & 0      & \cdots & 0      \\
\mb{A}_{2,1} & \mb{A}_{2,2} & 0      & \cdots & 0      \\
0      & \mb{A}_{3,2} & \mb{A}_{3,3} & \cdots & 0      \\
\vdots & \vdots & \vdots & \ddots & 0      \\
0      & 0      & 0      & \cdots & \mb{A}_{n_{t},n_{t}}
\end{bmatrix}^\intercal
\begin{bmatrix}
\boldsymbol{\psi}_1 \\
\boldsymbol{\psi}_2 \\
\boldsymbol{\psi}_3 \\
\vdots\\
\boldsymbol{\psi}_{n_{t}} \\
\end{bmatrix}
=
\begin{bmatrix}
\f{\p f}{\p \mb{x}_{1}} \\
\f{\p f}{\p \mb{x}_{2}} \\
\f{\p f}{\p \mb{x}_{3}} \\
\vdots \\
\f{\p f}{\p \mb{x}_{n_{t}}} \\
\end{bmatrix}
\renewcommand{\arraystretch}{1}
\end{equation}
where we can further decompose the equation with respect to ${\mb{x}}_{k}$ and ${\mb{P}}_{k}$:
\begin{align}
\mb{A}_{k, k} &=
\begin{bmatrix}
\label{eq:Akk}
\mb{A}_{k, k, \mb{x}, \mb{x}} & \mb{A}_{k, k, \mb{x}, \mb{P}} \\
\mb{A}_{k, k, \mb{P}, \mb{x}} & \mb{A}_{k, k, \mb{P}, \mb{P}}
\end{bmatrix}, \\[1em]
\mb{A}_{k, k-1} &=
\begin{bmatrix}
\label{eq:Akk-1}
\mb{A}_{k, k-1, \mb{x}, \mb{x}} & \mb{A}_{k, k-1, \mb{x}, \mb{P}} \\
\mb{A}_{k, k-1, \mb{P}, \mb{x}} & \mb{A}_{k, k-1, \mb{P}, \mb{P}}
\end{bmatrix}, \\[1em]
\boldsymbol{\psi}_k &=
\begin{bmatrix}
\boldsymbol{\psi}_{k, \mb{x}_k} \\
\boldsymbol{\psi}_{k, \mb{P}_k}
\end{bmatrix},  \label{eq:psi} \\[1em]
\frac{\partial f}{\partial \mb{x}_k} &=
\begin{bmatrix}
\frac{\partial f}{\partial \mb{x}_{k, \mb{x}_k}} \\
\frac{\partial f}{\partial \mb{x}_{k, \mb{P}_k}}
\end{bmatrix},  \label{eq:f}
\end{align}
where we can further decompose the equation with respect to the state $\mb{x}_k$ and the covariance $\mb{P}_k$. 
For convenience, we define the augmented variable
\begin{equation}
\label{eq:yk_def}
\mb{y}_k \coloneqq 
\begin{bmatrix}
\mb{x}_k \\
\vecop(\mb{P}_k)
\end{bmatrix}
\in \mathbb{R}^{\,n_x+n_p},
\end{equation}
where $\vecop(\cdot)$ denotes the column-wise vectorization of a matrix, 
$\mb{x}_k \in \mathbb{R}^{n_x}$ is the state vector, 
$\mb{P}_k \in \mathbb{R}^{n_x \times n_x}$ is the covariance matrix, 
and $n_p=n_x^2$ if the full matrix is vectorized.
Accordingly, the adjoint vector and the gradient of the objective with respect to $\mb{y}_k$ can be written in block form as
\begin{align}
\boldsymbol{\psi}_k &=
\begin{bmatrix}
\boldsymbol{\psi}_{k,\mb{x}} \\
\boldsymbol{\psi}_{k,\mb{P}}
\end{bmatrix}, \\[1em]
\frac{\partial f}{\partial \mb{y}_k} &=
\begin{bmatrix}
\dfrac{\partial f}{\partial \mb{x}_k} \\
\vecop\!\left(\dfrac{\partial f}{\partial \mb{P}_k}\right)
\end{bmatrix}, 
\end{align}
here, $\boldsymbol{\psi}_{k,\mb{x}} \in \mathbb{R}^{n_x}$ and ${\partial f}/{\partial \mb{x}_k} \in \mathbb{R}^{n_x}$ correspond to the state block, 
while $\boldsymbol{\psi}_{k,\mb{P}} \in \mathbb{R}^{n_p}$ and $\vecop({\partial f}/{\partial \mb{P}_k}) \in \mathbb{R}^{n_p}$ correspond to the covariance block.

The Eqs.~\eqref{eq:psi} and~\eqref{eq:f} are usually problem specific and can be easily computed.
For the block matrices in Eqs.~\eqref{eq:Akk} and \eqref{eq:Akk-1}, we have derived the analytic forms of them.
We present them here
\begin{equation}
\label{eq:sensitivity}
\begin{aligned}
\mb{A}_{k, k-1, \mb{x}, \mb{x}} &= - \left( \mathbf{I} - \mathbf{K}_{k} \mb{H}_{k} \right) \mathbf{F}_{k} \\
\mb{A}_{k, k-1, \mb{x}, \mathbf{P}} &= - \left[ \mb{S}_{k}^{-\intercal} \mb{H}_{k} \otimes (\mathbf{I} - \mathbf{K}_{k} \mb{H}_{k}) \right] 
(\mathbf{F}_{k} \otimes \mathbf{F}_{k}) (\mb{z}_{k} - \mb{H}_{k} \mb{x}_{k-1}) \\
\mb{A}_{k, k-1, \mathbf{P}, \mb{x}} &= \mb{0} \\
\mb{A}_{k, k-1, \mathbf{P}, \mathbf{P}} &= \left[ ( \mathbf{P}_{k-1}^\intercal \mb{H}_{k}^\intercal \mb{S}_{k}^{-\intercal} \mb{H}_{k} - \mathbf{I} ) \mathbf{F}_{k} \right]
\otimes \left[ (\mathbf{I} - \mathbf{K}_{k} \mb{H}_{k}) \mathbf{F}_{k} \right] \\
\mb{A}_{k, k, \mb{x}, \mb{x}} &= \mathbf{I} \\
\mb{A}_{k, k, \mb{x}, \mathbf{P}} &= \mb{0} \\
\mb{A}_{k, k, \mathbf{P}, \mb{x}} &= \mb{0} \\
\mb{A}_{k, k, \mathbf{P}, \mathbf{P}} &= \mathbf{I}
\end{aligned}
\end{equation}
We have verified these formulas numerically and the detailed derivation can be found in~\ref{sec:Sensitivity}. 

In summary, sensitivity computation using the adjoint method provides an efficient and robust means to evaluate the gradients needed in the optimization of the transition matrices $\hat{\mathbf{F}}_k$. 
Compared with automatic differentiation frameworks such as JAX~\cite{Jax2018a}, the analytic adjoint formulation avoids the overhead of unrolling the entire filtering process, thereby reducing memory cost and improving numerical stability, especially for long horizons. 
By combining these analytic sensitivities with the DKF structure, we enable scalable field inversion of dynamics. 
Next we will describe the closure model using machine learning in the following~\Cref{subsec:Closure model}.


\subsection{Closure model using machine learning}
\label{subsec:Closure model}
Once the system dynamics have been identified through field inversion in the first stage, we introduce a machine learning model to generalize these dynamics. 
Specifically, the transition operator is parameterized by a deep neural network
\begin{equation}
\hat{\mathbf{F}}_{\boldsymbol{\theta}} = \texttt{DNN}(\mathbf{x}, \mathbf{d}; \boldsymbol{\theta}),
\end{equation}
where $\boldsymbol{\theta}$ denotes the trainable weights and biases. 
The transition operator defines the model obtained through field inversion of the imperfect model in Eq.~\eqref{eq:imperfect_model}.
The network is trained to reproduce the optimized operators obtained in the first stage by minimizing
\begin{equation}
\begin{aligned}
\min \,\, & \text{loss}(\hat{\mathbf{F}}_{\boldsymbol{\theta}}, \hat{\mathbf{F}}) \\
\text{w.r.t.}\,\, & \boldsymbol{\theta},
\end{aligned}
\end{equation}
where the loss function is defined as the Frobenius norm of the difference,
\[
\texttt{loss}(\hat{\mathbf{F}}_{\boldsymbol{\theta}}, \hat{\mathbf{F}}) 
= \frac{1}{n_t} \sum_{k=1}^{n_t} 
\left\| \hat{\mathbf{F}}_{\boldsymbol{\theta}}(\mathbf{x}_k,\mathbf{d}) - \hat{\mathbf{F}}_k \right\|_F,
\]
with $n_t$ the number of time steps. 

Once this step is done, the discovered dynamics can be used to replace the original dynamics $\hat{\mb{F}}$ in the Kalman filter.
Here, the discovered dynamics are parameterized by a multilayer perceptron (MLP).

This two-stage workflow enables the integration of data-driven models into the classic Kalman filter framework. 
In the first phase, the underlying system dynamics are revealed through field inversion or direct optimization, providing an interpretable sequence of transition operators that best match the observed data. 
In the second phase, we leverage the expressive power of neural networks to learn a parameterized mapping from the latent state and design variables to the optimal transition matrices, effectively encoding complex, possibly nonlinear dependencies that are difficult to capture analytically.

The neural network $\texttt{DNN}(\mb{x}, \mb{d}; \boldsymbol{\theta})$ acts as a flexible surrogate for the transition operator, allowing for rapid, state-dependent updates during filtering. 
The training objective $\texttt{loss}(\hat{\mb{F}}_{\boldsymbol{\theta}}, \hat{\mb{F}})$ ensures that the learned model closely approximates the optimal operators identified in the first stage, while also enabling generalization to new states or operating regimes.

By embedding the discovered dynamics into the Kalman filter, we obtain a hybrid framework that combines the statistical rigor of state estimation with the adaptability of machine learning. 
This approach improves predictive accuracy, enables real-time adaptation to changing system conditions, and offers a scalable path for assimilating large, heterogeneous datasets. 
The effectiveness of this strategy will be further demonstrated in the following numerical experiments.

\section{Numerical experiments}
\label{sec:experiments}
To demonstrate the interpretability of the proposed DKF methodology, we conduct two numerical experiments, which include a classic rocket dynamics system in~\Cref{subsec:Rocket} and a nonlinear reaction-diffusion equation (the Allen--Cahn boundary value problem) in~\Cref{subsec:Allen--Cahn}. 
These examples are chosen to cover both low and high dimensional dynamical systems with different types of parameters and noise structures. 
The results highlight the robustness, accuracy, and flexibility of our approach for field inversion and dynamics discovery in complex, noisy environments. 
Detailed descriptions of the models, experimental setups, and comparative analyses are provided in the following subsections. 

\subsection{Rocket model}
\label{subsec:Rocket}
A simple test case for the DKF is a rocket launching problem.
The rocket has a thrust till time exceeds the burn time, then it behaves like a free-falling object under gravity. 
The governing equation is defined by
\begin{equation}
\label{eq:rocket_dynamics}
\hat{\mb{x}}_{k+1} = \hat{\mathbf{F}}_k \hat{\mb{x}}_k + \hat{\mb{B}}_k \hat{{u}}_k + \hat{\mb{f}}_k + \hat{\mathbf{w}}_k,
\end{equation}
where
\begin{align}
\hat{\mb{x}}_k &= 
\begin{bmatrix} x_k \\ \dot{x}_k \end{bmatrix}, \quad
\hat{\mathbf{F}}_k = 
\begin{bmatrix} 1 & \Delta t \\ 0 & 1 \end{bmatrix}, \quad
\hat{\mb{B}}_k = 
\begin{bmatrix} 0 \\ \frac{\Delta t F_{T}}{m} \end{bmatrix}, \\
\hat{\mb{f}}_k &= 
\begin{bmatrix} 0 \\ -g \Delta t \end{bmatrix}, \quad
\hat{{u}}_k = 
\begin{cases}
1, & \text{if } t_k < \text{burn time} \\
0, & \text{if } t_k > \text{burn time}
\end{cases}
\end{align}
$F_{T}$ is the thrust and $m$ is the mass,
and $\hat{\mathbf{w}}_k\in\mathbb{R}^{n_{{x}}}\sim \mathcal{N}(\mb{0}, \hat{\mb{Q}}_k)$ denotes the external noise, $\hat{\mb{Q}}_k$ is the covariance matrix of the observation noise.
We measure the altitude of the rocket, so the observation operator becomes:
\begin{equation}
\hat{\mb{H}}_k = \begin{bmatrix} 1 \quad 0 \end{bmatrix},
\end{equation}
Here, the transition matrix $\hat{\mathbf{F}}_k$ is known. 
In a general case, this matrix may not be known with certainty, that is, cases with unknown or partially known dynamics. 
In such a case, we can treat the matrices $\hat{\mathbf{F}}_k$ as unknown design variables. 
And the optimization problem as described in Eq.~\eqref{eq:opt} becomes: find the matrices ${\mathbf{F}}_k$, across all time steps, such that the estimated value of $\hat{\mb{x}}_k$ becomes close to the measured value.



We consider a dynamical system governed by physical processes such as thrust and gravity, where the goal is to estimate the system state from noisy observations. 
In this setting, the state transition matrix $\hat{\mathbf{F}}_k$ at each time step is assumed unknown and is treated as a design variable to be identified. 
The observation operator and noise covariance matrices are specified according to the system's physical characteristics.

Given a sequence of noisy observation $\{\hat{\mb{z}}_k\}_{k=1}^{n_{t}}$, the objective is to determine the transition matrices $\hat{\mathbf{F}}_k$ that minimizes the discrepancy between the predicted (filtered) states and the observations over the trajectory. 
The estimation problem is formulated as the following optimization:
\begin{equation}
\begin{aligned}
\min \,\, & \sum_{k=1}^{n_{t}} \left\| \hat{\mb{H}}_{k} \hat{\mb{x}}_{k}- \hat{\mb{z}}_{k} \right\|^2 \\
\text{w.r.t.}\,\, & \hat{\mathbf{F}}_1, \ldots, \hat{\mathbf{F}}_{n_t},
\end{aligned}
\end{equation}
where $\hat{\mb{x}}_{k}$ denotes the posterior state estimate at each filtering step, $\hat{\mb{z}}_{k}$ is the observed data, and $\hat{\mb{H}}_{k}$ is the observation matrix. 
The full parameter vector is constructed by concatenating the sequence of $\hat{\mathbf{F}}_k$ matrices across all time steps.
We will solve this optimization problem in the next subsection.

\subsubsection{Solution of optimization problem}

To efficiently solve this optimization problem, we leverage the DKF framework and compute gradients of the loss function with respect to $\mb{d}$ using the adjoint method. 
Analytic Jacobian block matrices are derived and assembled based on theoretical results to ensure numerical stability and accurate gradient evaluation. 
The sequence ${\hat{\mathbf{F}}_k}$ is initialized as small perturbations around a physically reasonable baseline, and iterative optimization is performed using the L-BFGS algorithm.

~\Cref{tab:F_compare} presents a comparison of the state transition matrices at selected time steps, showing the true values, initial values, and the optimized results. 
The optimized $\mathbf{F}_k$ matrices consistently approach the true dynamics, demonstrating the effectiveness of the inversion process.

\begin{figure}[htbp]
\centering
\includegraphics[width=\textwidth]{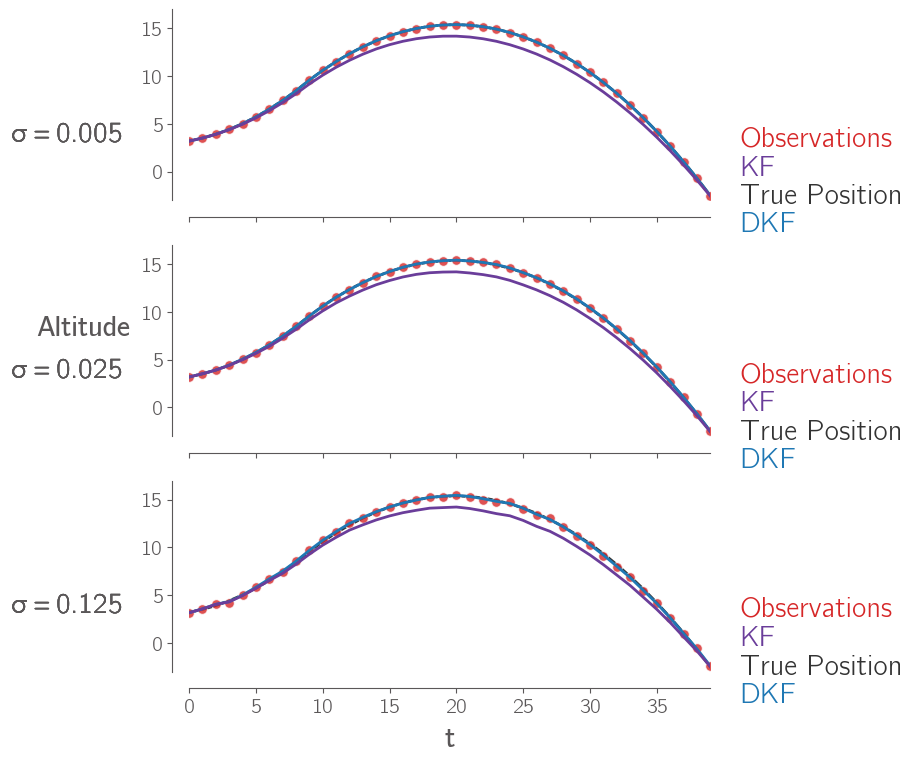}

\caption{Optimization results for three levels of observation noise: $\sigma = 0.005, 0.025, 0.125$.}
\label{fig:rocket_noise_results}
\end{figure}

Finally, the impact of the optimized parameters on filtering performance is systematically evaluated under various observation noise levels.
The black solid line (``True Position'') shows the ground truth trajectory of the system generated by the dynamics. 
The red dots (``Observations'') indicate the observed values at each time step, which are influenced by Gaussian noise with the standard deviation $\sigma$ specified in each subplot. 
The blue line (``DKF'') corresponds to the filtered state estimates obtained by the DKF using the optimized transition matrix sequence $\mathbf{F}_k$. 
The purple solid line (``KF'') represents the filter output when using the initial $\hat{\mathbf{F}}_k$ sequence without optimization. 
The results in Fig.~\ref{fig:rocket_noise_results} illustrate the filtering performance of the optimized system under various observation noise levels.

\begin{table}[htbp]
\centering
\caption{Comparison of the global state transition matrix $\mathbf{F}_k$ (true, initial, and optimized) for different observation noise levels.}
\label{tab:F_compare}
\resizebox{\textwidth}{!}{
\begin{tabular}{lccc}
\toprule
$\sigma$ & True $\mathbf{F}_k$ & Initial $\hat{\mathbf{F}}_k$ & Optimized $\mathbf{F}_k$ \\
\midrule
$0.005$ &
$\begin{bmatrix}
1.000000 & 0.100000 \\
0.000000 & 1.000000
\end{bmatrix}$ &
$\begin{bmatrix}
0.957691 & 0.088596 \\
-0.072960 & 0.967528
\end{bmatrix}$ &
$\begin{bmatrix}
1.002941 & 0.100005 \\
-0.000087 & 0.997059
\end{bmatrix}$ \\
\midrule
$0.025$ &
$\begin{bmatrix}
1.000000 & 0.100000 \\
0.000000 & 1.000000
\end{bmatrix}$ &
$\begin{bmatrix}
0.957691 & 0.088596 \\
-0.072960 & 0.967528
\end{bmatrix}$ &
$\begin{bmatrix}
1.003099 & 0.097841 \\
-0.000109 & 0.997220
\end{bmatrix}$ \\
\midrule
$0.125$ &
$\begin{bmatrix}
1.000000 & 0.100000 \\
0.000000 & 1.000000
\end{bmatrix}$ &
$\begin{bmatrix}
0.957691 & 0.088596 \\
-0.072960 & 0.967528
\end{bmatrix}$ &
$\begin{bmatrix}
1.003100 & 0.097842 \\
-0.000109 & 0.997221
\end{bmatrix}$ \\
\bottomrule
\end{tabular}}
\end{table}

As the noise standard deviation increases, the tracking accuracy of the filter generally decreases. 
However, for all noise levels, the optimized $\mathbf{F}_k$ sequence consistently provides a closer fit to the true state than the unoptimized $\mathbf{F}_k$ sequence, demonstrating the effectiveness of the optimization procedure. 
In particular, at low noise levels ($\sigma=0.005$ and $\sigma=0.025$), the optimized filter nearly overlaps with the true state trajectory, whereas the initial estimate shows noticeable bias. 
Even as the noise increases to $\sigma=0.125$, the optimized filter maintains a robust performance and clearly outperforms the unoptimized case. 

While the current study focuses on a low-dimensional rocket model, in practice, many dynamical systems of interest are high-dimensional, resulting in a large covariance matrix $\mathbf{P}_k$.
For high-dimensional systems, direct update of the full covariance matrix $\mathbf{P}_k$ is computationally expensive and often infeasible. 
To overcome this, one effective approach is to use a block-diagonal or localized structure for $\mathbf{P}_k$, which only tracks correlations between closely related variables and ignores weak or distant dependencies, thereby significantly reducing computational load. 
Another approach is to employ low-rank or subspace approximations—such as principal component analysis (PCA) or singular value decomposition (SVD)—to represent $\mathbf{P}_k$ in terms of its most important modes. 

These results confirm that the proposed method can effectively learn the underlying system parameters and improve filtering accuracy, even in the presence of significant observation noise.

\subsection{Allen--Cahn boundary value problem}
\label{subsec:Allen--Cahn}
The Allen--Cahn equation is a prototypical nonlinear reaction–diffusion model that captures phase separation and interface dynamics in multi-component systems \cite{Zhang2009a,Qi2023a,Mohammadi2018a,Benner2013a}. 
In this section, we consider a boundary value problem governed by a generalized Allen--Cahn equation with nonlinear, state-dependent diffusion. 
Our aim is to recover the underlying dynamics and explore robust state estimation under various observation noise levels.

This section is organized as follows. 
We first present the mathematical formulation of the governing Allen--Cahn equation in~\Cref{subsubsec:Govern}. 
We then describe the solution scheme of the boundary value problem in~\Cref{subsubsec:Solution method}. 
A central focus is placed on the DKF, which enables direct learning of the system's underlying dynamics from noisy observations~\Cref{subsubsec:Learning the underlying dynamics via DKF}. 
Finally, we compare the prediction accuracy by examining both the mean and variance of their estimated states in~\Cref{subsubsec:Prediction of mean} and~\Cref{subsubsec:Prediction of variance}. 

\subsubsection{Governing equation}
\label{subsubsec:Govern}
The system is governed by a nonlinear reaction–diffusion equation with state-dependent diffusivity:
\begin{equation}
\frac{\partial v}{\partial t} = \frac{\partial}{\partial x} \left[ d(v) \frac{\partial v}{\partial x} \right] - v^3 + v, \quad x \in [0, l],\ t > 0,
\end{equation}
where $d(v) = 0.1\, \tanh(v)$. We impose periodic boundary conditions to ensure both the solution and its derivative are continuous at the domain endpoints:
\begin{equation}
v(0, t) = v(l, t), \quad \frac{\partial v}{\partial x}\Big|_{x=0} = \frac{\partial v}{\partial x}\Big|_{x=l}, \quad t > 0.
\end{equation}
The initial condition is a stripe pattern with a small sinusoidal perturbation, projected onto $[0,1]$:
\begin{align}
b(x) &\coloneqq \frac{1+\operatorname{sgn}\!\big(\frac{\sin(2\pi m_{\mathrm{stripes}})x}{l}\big)}{2}\in\{0,1\},\\
\tilde {v}_0(x) &\coloneqq b(x) + A\,\sin\!\Big(2\pi k\,\frac{x}{l} + \phi\Big),\\
v_0(x) &\coloneqq \Pi_{[0,1]}\big(\tilde {v}_0(x)\big),
\end{align}
where $A=0.18$, $k=3$, $\phi\sim\mathcal U(0,2\pi)$, 
$m_{\mathrm{stripes}}\in\mathbb N$ (we use $m_{\mathrm{stripes}}=8$), and 
$\Pi_{[0,1]}(w)\coloneqq\min\{1,\max\{0,w\}\}$ denotes projection onto $[0,1]$ (here $w$ is only a dummy variable used to define the projection operator; in practice it is replaced by the specific argument, e.g.\ $\tilde v_0(x)$).
Here, $b(x)$ is a $0/1$ stripe template marking the two phases; $\tilde {v}_0(x)$ adds a small within-stripe sinusoidal variation; and $v_0(x)$ is the clipped initial condition enforcing bounds $[0,1]$.

\subsubsection{Solution method}
\label{subsubsec:Solution method}
This experiment focuses on~\Cref{subsubsec:Govern}. 
The system is discretized over $N=16$ spatial grid points and $T=30$ time steps. 
The spatial mesh size is $dx = l/N$, and the time step is $dt = 0.01$. 

To investigate this inverse problem, we generate synthetic data by time marching the conservative form of the reaction–diffusion dynamics with periodic boundaries. At each time step we first evaluate the nonlinear diffusivity
\(d_i^n \coloneqq d(v_i^n) = 0.1\,\tanh(v_i^n)\) at cell centers and then form face coefficients by arithmetic averaging
\(d_{i+\frac{1}{2}}^n \coloneqq {1}/{2}\,(d_i^n + d_{i+1}^n)\) (periodic indexing).
The state is advanced by an explicit Euler scheme using fluxes
\begin{equation}
\begin{aligned}
J_{i+\frac{1}{2}}^n &= d_{i+\frac{1}{2}}^n 
   \frac{(v_{i+1}^n - v_i^n)}{\Delta x},\\
v_{i}^{n+1}
&= v_i^n + \Delta t \left[
\frac{J_{i+\frac{1}{2}}^n - J_{i-\frac{1}{2}}^n}{\Delta x}
- (v_i^n)^3 + v_i^n
\right] \\
&= v_i^n + \Delta t \left[
\frac{ d_{i+\frac{1}{2}}^n (v_{i+1}^n - v_i^n) - d_{i-\frac{1}{2}}^n (v_i^n - v_{i-1}^n)}{(\Delta x)^2}
- (v_i^n)^3 + v_i^n
\right],
\end{aligned}
\end{equation}
gaussian noise with varying standard deviation is added independently to all spatial points at each time step to generate the observation sequence $\hat{\mb{z}}_t$.

A key feature of our approach is that, instead of assuming the system dynamics are fully known, we leverage the DKF to simultaneously perform filtering and operator learning. 
The diffusivity $d(v)$ is represented as independent values on a uniform grid of state amplitudes, allowing local adaptation to different regions of the state space. 
Specifically, we introduce tabulated parameters $\tilde{d}(v_j)$, where each $\tilde{d}(v_j)$ denotes the learned diffusivity value associated with grid point $v_j$ in the state domain. 
These table values $\{\tilde{d}(v_j)\}$ are optimized to minimize the mean squared innovation between predicted and observed data.
Formally, the loss is defined as
\begin{equation}
\begin{aligned}
\min \,\, & \frac{1}{n_{t}} \sum_{t=1}^{n_{t}} 
\big\| \mb{S}_t^{-\frac{1}{2}}\big(\mb{z}_t - \mb{H} x_t^{\mathrm{}}\big) \big\|^2 \\
\text{w.r.t.}\,\, & \tilde{d}(v_1),...,\tilde{d}(v_{n_{t}}), 
\end{aligned}
\end{equation}
where $\mb{S}_t$ is the innovation covariance from the Kalman filter. 
Whitening by $\mb{S}_t^{-1/2}$ ensures statistical consistency by accounting for the varying uncertainty at each time step. 
The optimized diffusivity table then defines the operator sequence $\mathbf{F}_k$, which provides a data-driven approximation of the true dynamics and is used to reconstruct the state trajectory.

In addition, we train a neural network surrogate to predict $d(v)$ directly from the state $v$. 
The network is trained on the field-inversion results by minimizing the mean squared error between its predicted diffusivity and the reference table values (with preprocessing to improve robustness).
Once trained, the neural network–based $d(v)$ is used to reconstruct the state evolution over the entire time window.








\subsubsection{Learning the underlying dynamics via DKF}
\label{subsubsec:Learning the underlying dynamics via DKF}

A central aim of this work is to learn the underlying dynamics of the system directly from noisy observations. Here, ``dynamics'' refer to the latent operator $\hat{\mathbf{F}}_k$ that governs the system’s evolution at each time step.
While the classic Kalman filter typically assumes these dynamics are known and fixed, our DKF framework enables adaptive, data-driven discovery of the governing operator $\hat{\mathbf{F}}_k$. 

More specifically, instead of directly optimizing the full discrete evolution operator $\hat{\mathbf{F}}_k$, we adopt a physics-informed approach by optimizing the nonlinear diffusivity $d(v)$—using either an explicit parametric form, a nonparametric point-wise representation on a uniform $u$-grid, or, in post-processing, a neural network fitted to the inferred $d(v)$.
This ensures that the learned dynamics remain physically consistent (\emph{e.g.}, $d(v) > 0$) and interpretable.
The system operator $\hat{\mathbf{F}}_k$ is then constructed as a function of $d(v)$, thus ensuring that the learned dynamics remain physically consistent and interpretable. 
This strategy allows DKF to actively learn and refine the true physical law driving the process, even in the presence of substantial observational noise.

At each time step $k$, the operator $\hat{\mathbf{F}}_k$ is constructed as
\begin{equation}
\hat{\mathbf{F}}_k = \mathbf{I} + \Delta t \left[ \mathrm{diag}(d({v}_k))\, \frac{{d}_2}{\Delta x^2} + \mathrm{diag}(1 - 3v_k^2) \right],
\end{equation}
optionally with an additive bias term 
\(
b_k = 2 \Delta t\, v_k^3
\)
arising from the nonlinear reaction term. 
Here, $\mathbf{I}$ is the identity matrix, and ${d}_2$ is the central-difference Laplacian with periodic boundary conditions:
\begin{equation}
({d}_2)_{ij} =
\begin{cases}
1,  & j = i-1 \mod N, \\
-2, & j = i, \\
1,  & j = i+1 \mod N, \\
0,  & \text{otherwise},
\end{cases}
\end{equation}
where $N$ is the number of spatial grid points. 
This structure ensures that the estimated dynamics remain consistent with the underlying physical process, while reducing the dimensionality and improving the interpretability of the inverse problem.

\begin{figure}[htbp]
\centering
\includegraphics[width=1.\textwidth]{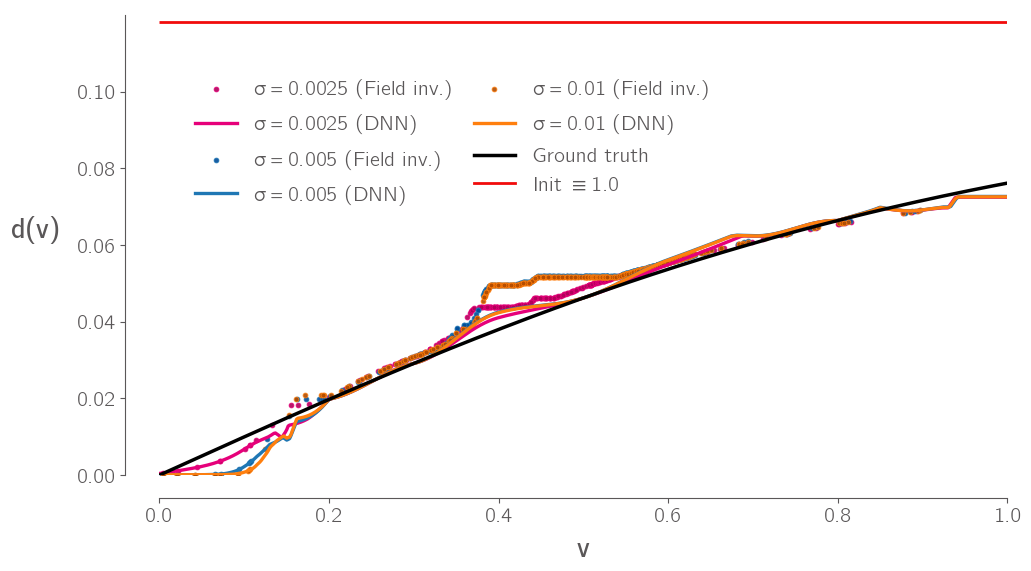}
\caption{
Comparison of $d(v)$ estimated by different methods under different noise levels $\sigma = 0.0025$, $0.005$, and $0.01$.
Dots indicate field inversion results, solid lines indicate DNN fits, 
the black line shows the ground truth $d(v)=0.1\tanh(v)$, 
and the red line denotes the constant initial guess $d(v)\equiv 1.0$. 
}
\label{fig:diffusion_compare}
\end{figure}

The results in Fig.~\ref{fig:diffusion_compare} highlight the effectiveness of our framework in recovering the underlying diffusion law. 
First, the pointwise inversion performed via DKF optimization successfully reconstructs the nonlinear diffusivity $d(v)$ across all noise levels. 
Even under relatively large observational noise, the inverted values preserve the correct functional form and monotonic trend of $d(v)$, closely following the ground truth. 
This confirms that the DKF can directly infer reliable operator values from noisy and partial measurements without relying on a restrictive parametric assumption.  

On top of the inversion, we further train deep neural networks to approximate $d(v)$ using the DKF-inferred trajectories. 
The trained DNN models not only reproduce the inversion results with high fidelity but also yield smooth, continuous reconstructions that enhance generalization beyond the sampled states. 
Importantly, the DNN outputs remain consistent with the physical constraints (e.g., $d(v) > 0$) while providing a compact and interpretable representation of the learned dynamics.  

Together, these results demonstrate that our two-stage strategy---pointwise DKF inversion followed by DNN training---achieves both accuracy and robustness: the inversion ensures faithful recovery of the true operator from noisy data, while the DNN serves as a stable surrogate that captures the same physics in a smooth and generalizable manner.

To quantitatively evaluate the quality of this learned operator, we compare the inferred $\hat{\mathbf{F}}_k$ at each time step against the true operator $\mathbf{F}_k^{\mathrm{true}}$ used to generate the synthetic data. 
Specifically, we report the relative Frobenius norm error:
\[
\frac{\|\hat{\mathbf{F}}_k - \mathbf{F}_k^{\mathrm{true}}\|_F}{\|\mathbf{F}_k^{\mathrm{true}}\|_F}
\]
as a function of time, for both the initial guess and the optimized result after field inversion and DNN training.

\begin{figure}[htbp]
\centering
\includegraphics[width=1.\textwidth]{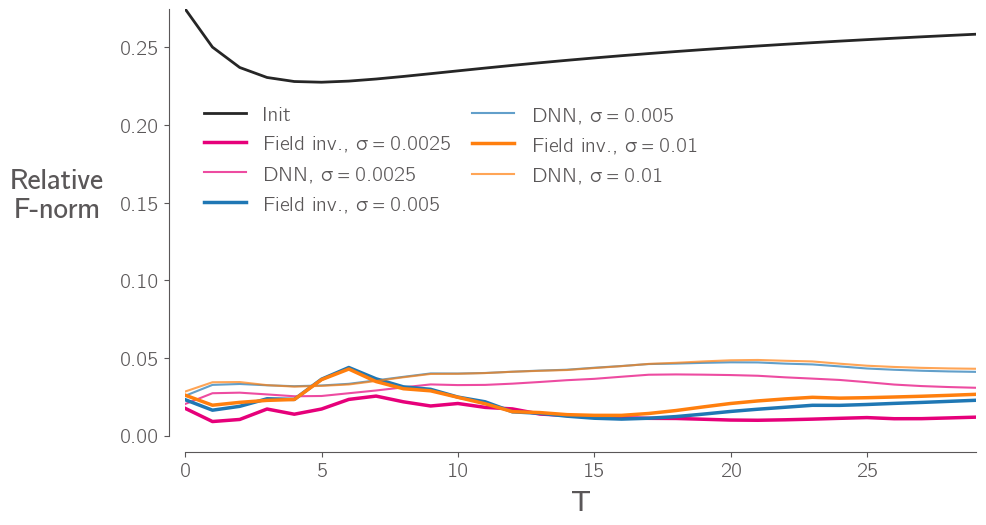} 

\caption{
Time evolution of the relative Frobenius norm error
${\|\hat{\mathbf{F}}_k - \mathbf{F}_{k}^{\mathrm{true}}\|_F}/{\|\mathbf{F}_k^{\mathrm{true}}\|_F}$
between the learned dynamical operator and the ground truth, for three levels of observation noise: 
$\sigma = 0.0025$, $0.005$, and $0.01$. 
}
\label{fig:F_norm_error_combined}
\end{figure}

As shown in Fig.~\ref{fig:F_norm_error_combined}, 
the initial transition matrix (black line) exhibits a large and persistent error with respect to the true dynamics, remaining around $0.25$ over the entire time horizon. 
By contrast, both the inversion results (solid colored lines) 
and the DNN predictions (lighter lines) achieve much smaller errors, consistently staying within the range $0.00$--$0.05$. 
This demonstrates that the proposed framework consistently learns transition matrix that accurately approximates the true dynamics across all tested noise levels.
Nevertheless, the DKF framework achieves a substantial error reduction compared to the initial guess, even at the highest noise level.

These results demonstrate that the DKF approach not only enables effective filtering and smoothing, but also robust, data-driven identification of underlying dynamics from noisy measurements. 
The learned time-varying transition matrix $\mathbf{F}_k$ provide interpretable, physically meaningful surrogates for the unknown evolution law, even when direct physical modeling is unavailable or inaccurate.

\subsubsection{Prediction of mean}
\label{subsubsec:Prediction of mean}
\begin{figure}[htbp]
\centering
\includegraphics[width=\textwidth]{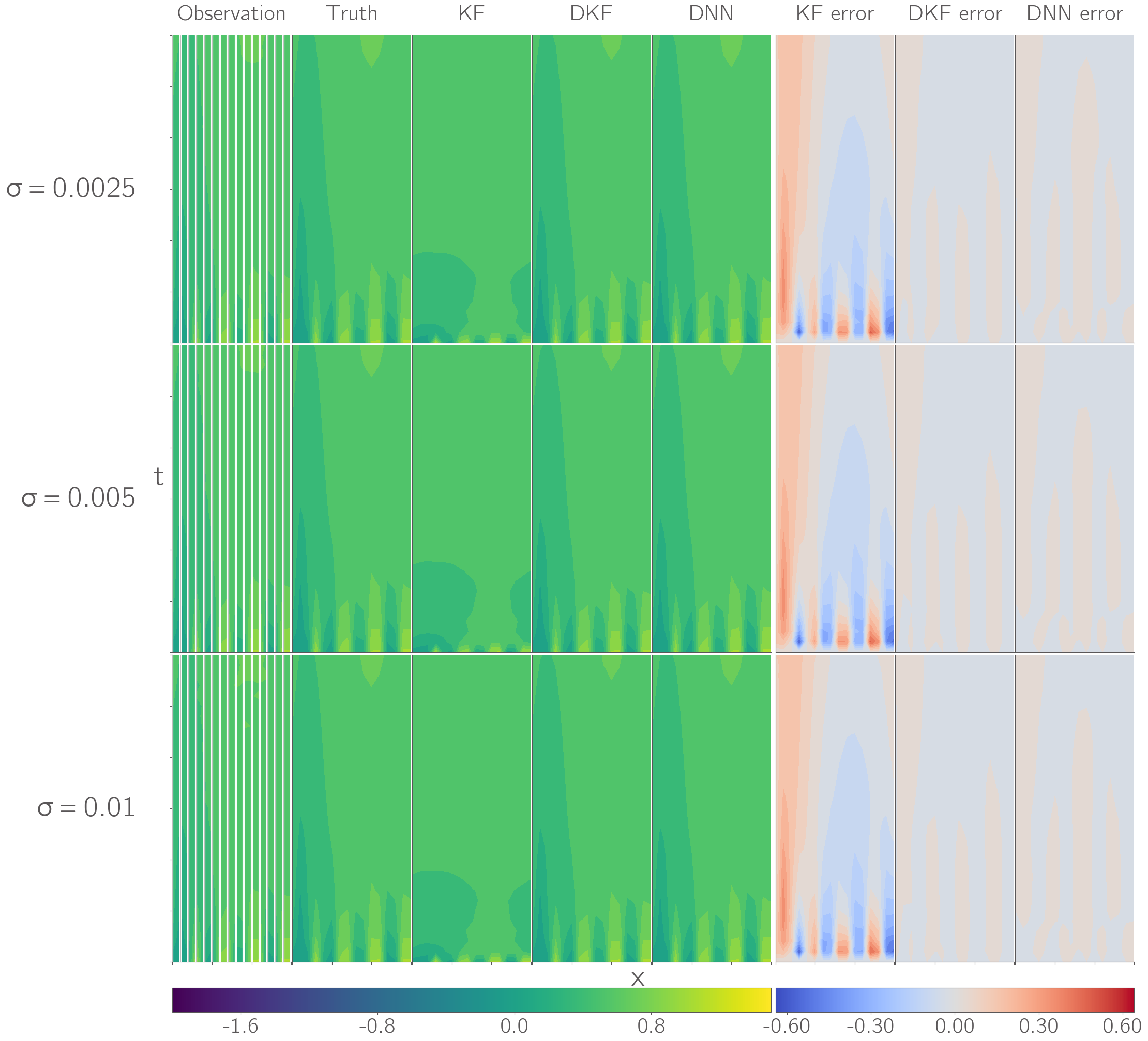}
\caption{
Spatiotemporal reconstruction under different observation noise levels 
($\sigma=0.0025$, $0.005$, $0.01$; rows). 
Columns show sparse observations $\hat{\mb{z}}_k$, ground truth $v_{\text{true}}^k$, 
Kalman filter without inversion $v_{\text{no\_inv}}^k$, DKF inversion $v_{\text{inv}}^k$, 
DNN predictions $v_{\text{DNN}}^k$, and the corresponding error fields. 
White vertical lines in the first column mark observed locations.
}
\label{fig:allen_cahn_subset_results}
\end{figure}

To evaluate the benefits of field inversion and machine learning in mean state estimation, we compare spatiotemporal reconstructions under varying observation noise, as shown in Fig.~\ref{fig:allen_cahn_subset_results}. 
Each row in Fig.~\ref{fig:allen_cahn_subset_results} presents the results for a different noise level ($\sigma=0.0025$, $0.005$, $0.01$), with columns displaying, 
from left to right, the sparse observations $\hat{\mb{z}}_k$, the ground truth ${v}_{\text{true}}^k$, the classic Kalman filter predictions without inversion $v_{\text{no\_inv}}^k$, the DKF field inversion $v_{\text{inv}}^k$, the DNN predictions $v_{\text{DNN}}^k$, and the corresponding absolute error fields $|v_{\text{inv}}^k - v_{\text{true}}^k|$ and $|v_{\text{DNN}}^k - v_{\text{true}}^k|$. 
In the first column, the sparse observations are illustrated by overlaying vertical white lines on the spatiotemporal field, clearly indicating the specific grid points where measurements are available, while the rest of the domain remains unobserved.

As the noise increases, the classic Kalman filter, which relies on a fixed and potentially mismatched dynamical operator, gradually accumulates bias and exhibits amplified errors—particularly in regions of sharp gradients and in later time steps. 
This is visually apparent in both the reconstructed profiles and the error fields, where the classic Kalman filter struggles to track the true state and produces distorted mean profiles when the signal-to-noise ratio is low.


In contrast, the DKF is able to directly correct the underlying dynamics by learning from data, significantly reducing both bias and variance in the estimated states. 
The error fields confirm that the DKF approach not only suppresses noise amplification but also effectively eliminates large-scale model mismatch.
The advantage of this inversion step is particularly prominent at moderate and relatively high noise levels, where operator correction plays a crucial role in maintaining robust estimation accuracy.

Building upon the improved trajectories from field inversion, the subsequent machine learning step further enhances predictive performance. 
The neural network surrogate, trained on DKF-optimized samples, accurately recovers the underlying diffusion operator and achieves near-perfect agreement with the true mean state—even in the presence of sparse and noisy observations. 
This data-driven approach not only matches the field inversion result but also demonstrates strong generalization. As a result, the closure model closely follows the ground truth and outperforms all other methods under increasing noise.

\subsubsection{Prediction of variance}
\label{subsubsec:Prediction of variance}

\begin{figure}[htbp]
\centering
\includegraphics[width=\textwidth]{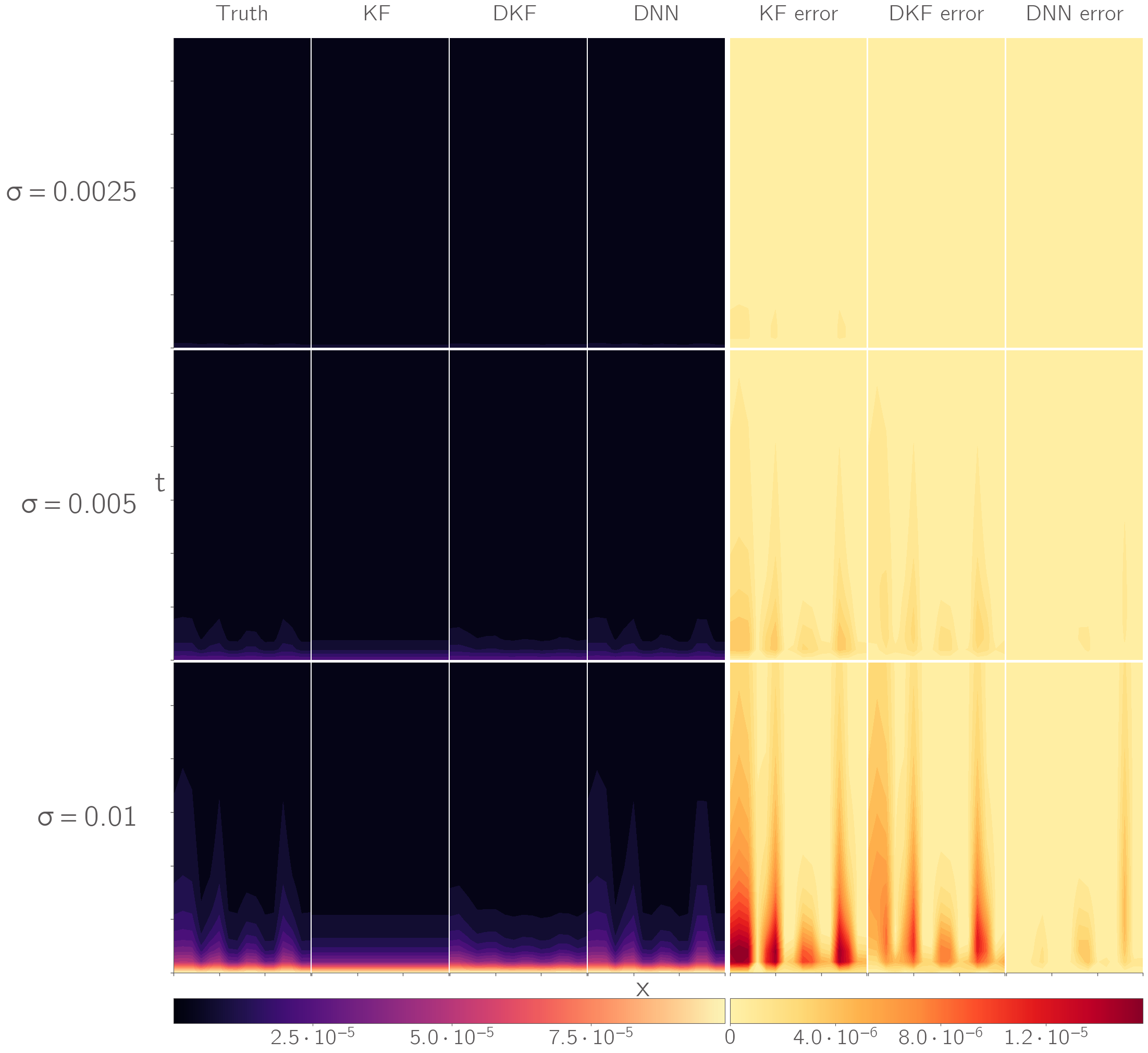}
\caption{
Evolution of the diagonal elements of the covariance matrix $\mathbf{P}_k$ 
and the corresponding absolute estimation errors under different observation noise levels 
($\sigma=0.0025$, $0.005$, $0.01$).
Panels show the ground truth, Kalman filter without inversion, DKF inversion, DNN prediction, and their respective errors. 
}
\label{fig:allen_cahn_subset_pk_results}
\end{figure}

Accurate estimation of the posterior covariance matrix is essential for quantifying uncertainty in nonlinear state estimation tasks.
The sequence of covariance matrices $\mathbf{P}_k$ provides not only a measure of confidence in the predicted states, but also enables principled uncertainty propagation for downstream inference or control.

In Fig.~\ref{fig:allen_cahn_subset_pk_results}, we present a comprehensive comparison of posterior variance estimation and error quantification under different observation noise levels ($\sigma=0.0025, 0.005, 0.01$). 
The first panel serves as the reference, revealing the true spatial and temporal structure of the variance field $[\mathbf{P}^k_{\text{true}}]_{ii}$. 
The second panel shows the posterior variance obtained from the classic Kalman filter without field inversion $[\mathbf{P}^k_{\text{no\_inv}}]_{ii}$, while the third panel displays the posterior variance recovered via differentiable field inversion $[\mathbf{P}^k_{\text{inv}}]_{ii}$. 
The fourth panel reports the posterior variance predicted directly by the deep neural network $[\mathbf{P}^k_{\text{DNN}}]_{ii}$. 
To facilitate direct comparison against the ground truth, the fifth panel visualizes the absolute error of the no-inversion estimate $\big|[\mathbf{P}^k_{\text{no\_inv}}]_{ii} - [\mathbf{P}^k_{\text{true}}]_{ii}\big|$, the sixth panel shows the error of the DKF estimate $\big|[\mathbf{P}^k_{\text{inv}}]_{ii} - [\mathbf{P}^k_{\text{true}}]_{ii}\big|$, and the seventh panel displays the corresponding error for the neural network prediction $\big|[\mathbf{P}^k_{\text{DNN}}]_{ii} - [\mathbf{P}^k_{\text{true}}]_{ii}\big|$. 
At the lowest noise level ($\sigma=0.0025$), all three estimation strategies produce variance maps that are nearly indistinguishable from the ground truth, as confirmed by the very low and spatially uniform absolute errors in the rightmost panels.
As the noise level increases, the limitations of the classic Kalman filter (no inversion) become apparent: its posterior variance exhibits systematic deviations from the true uncertainty, leading to persistent and sometimes spatially structured errors, especially in regions of high variability or low information.


In contrast, both the DKF and DNN approaches maintain high-fidelity variance estimation even under substantial noise. The DKF leverages data-driven operator correction to dynamically adapt its uncertainty quantification, closely matching the true variance in both magnitude and spatial pattern. 
The DNN surrogate, trained on DKF-optimized operator, inherits and even strengthens this capability, achieving uniformly low absolute errors across the entire domain. 
This robust performance highlights the critical advantage of field inversion and machine learning in capturing not only the mean but also the variance associated with each state estimate.

Overall, these results demonstrate that, while classic Kalman filtering may significantly misestimate uncertainty when the model is mismatched or noise is high, field inversion and machine learning can provide trustworthy, spatially resolved quantification of prediction confidence—a capability essential for reliable scientific inference and downstream decision-making in partially observed dynamical systems.

\section{Conclusion}
In this paper, we presented a differentiable Kalman filtering (DKF) framework that integrates physics-based state estimation with machine learning–driven model discovery. 
Classic Kalman filters are constrained by the assumption of fixed, often imperfect dynamics. 
By introducing a two-level adjoint-based optimization procedure, our approach—-DKF-based field inversion coupled with machine learning—adaptively corrected the dynamics model using observed data, while leveraging neural operators for enhanced predictive capability.

We demonstrated the effectiveness of this framework on two representative benchmarks: a rocket model with algebraic dynamics and a nonlinear reaction–diffusion PDE (the Allen–Cahn equation). 
In both cases, our method consistently outperformed the classic Kalman filter in terms of state reconstruction accuracy and robustness to observation noise. 
For example, in the Allen–Cahn system, we achieved a root mean square error (RMSE) of $d(v)$ below $10^{-2}$ across moderate noise levels ($\sigma = 0.0025$ to $0.01$), representing an improvement of at least two orders of magnitude compared to the classic Kalman filter. 
Furthermore, our deep neural network, trained on field-inverted trajectories, reduced estimation errors even further and exhibited strong generalization under substantial observational noise.

Beyond accuracy, the framework provided interpretable uncertainty quantification through principled covariance propagation. 
The field inversion step yielded physical insight into the governing dynamics, while the machine learning component enhanced flexibility and generalization.

Future work will extend this framework to more complex PDEs and multi-physics systems, incorporating Bayesian priors and ensemble-based approaches to further improve uncertainty quantification and enable real-time inference.


\bibliographystyle{elsarticle-num}
\bibliography{bib/mdolab,bib/references}

\setcounter{secnumdepth}{3}
\appendix
\section{Kalman filter derivation}
\label{sec:Kalman filter derivation}
In this section, we present the detailed derivation of the Kalman filter related formulas. 
We follow the equations and definitions as mentioned in Eq.~\eqref{eq:imperfect_model}.

Kalman filter operates with a given model of the underlying dynamics.
The Kalman filter operates in two steps: \textit{predict} and \textit{update}.
In the predict step, the filter estimates the current state $\hat{\mb{x}}_{k|k-1}$ and the error covariance $\mb{P}_{k|k-1}$ based on the previous  $\hat{\mb{x}}_{k-1}$ and the previous error covariance $\mb{P}_{k-1}$:
\begin{equation}
\label{eq:Kalman_predict}
\begin{aligned}
\hat{\mb{x}}_{k|k-1} &= \hat{\mb{F}}_k \hat{\mb{x}}_{k-1} + \hat{\mb{B}}_k \hat{\mb{u}}_k + \hat{\mb{f}}_k, \\
\mb{P}_{k|k-1} &= \hat{\mb{F}}_k \mb{P}_{k-1} \hat{\mb{F}}_k^\intercal + \hat{\mb{Q}}_k,
\end{aligned}
\end{equation}
where $\hat{\mb{x}}_{k|k-1}$ is the predicted state estimate at time step $k$, and $\mb{P}_{k|k-1}$ is the predicted error covariance.

In the update step, the filter incorporates the new observation $\hat{\mb{z}}_k$ to refine the state estimate $\hat{\mb{x}}_k$ and the error covariance $\mb{P}_k$:
\begin{equation}
\label{eq:Kalman_update}
\begin{aligned}
\hat{\mb{x}}_k &= \hat{\mb{x}}_{k|k-1} + \mb{K}_k \left(\hat{\mb{z}}_k - \hat{\mb{H}}_k \hat{\mb{x}}_{k|k-1}\right),\\
\mb{P}_k &= (\mb{I} - \mb{K}_k \hat{\mb{H}}_k) \mb{P}_{k|k-1},\\
\mb{K}_k &= \mb{P}_{k|k-1} \hat{\mb{H}}_k^\intercal \left(\hat{\mb{H}}_k \mb{P}_{k|k-1} \hat{\mb{H}}_k^\intercal + \hat{\mb{R}}_k \right)^{-1}, 
\end{aligned}
\end{equation}
where $\mb{K}_k$ is the Kalman gain, $\mb{I}$ is the identity matrix, $\hat{\mb{z}}_k$ is the observation at time step $k$, $\hat{\mb{x}}_k$ is the updated state estimate, and $\mb{P}_k$ is the updated error covariance.

The predict step Eq.~\eqref{eq:Kalman_update} can be written in their residual form by simply transferring all the RHS terms to LHS and equating it to zero that gives the residual Eq.~\eqref{eq:KF_residual}. 
To solve this problem, an initial guess for the solution ${\hat{\mb{x}}}_{k} \text{ and } \mb{P}_{k}$ can be taken and a newton method-based root finding algorithm can be used to converge to the actual solution. 
Therefore, at each step, ${\hat{\mb{x}}}_{k} \text{ and } \mb{P}_{k}$ becomes the state variables. 
In other words, the flattened matrices $[{\hat{\mb{x}}}_{1}, \mb{P}_{1}, {\hat{\mb{x}}}_{2}, \mb{P}_{2}, \ldots ]$ become the total state vector of this problem. 

We can take Kalman filter as a paramitrized mapping of $\left(\hat{\mb{x}}_0, \mb{P}_0, \hat{\mb{Z}}\right)\rightarrow \left(\hat{\mb{x}}_k, \mb{P}_k\right)$ where $\hat{\mb{Z}}=\begin{bmatrix}{\hat{\mb{z}}}_1, \ldots, {\hat{\mb{z}}}_k\end{bmatrix}$ is the observation matrix comes by observing the underlying physics.
To differentiate this equation using algorithmic differentiation tools, it is helpful to look at in intermediate steps and visualize how the intermediate variables depend on each other. 
Computing partial derivatives using reverse algorithmic differentiation for adjoint method involves seeding (or perturbing) the output (the residuals) and seeing how that effect propagates back to the inputs.

\section{Sensitivity computation}
\label{sec:Sensitivity}

According to Eq.~\eqref{eq:KF_residual}, $\mb{r}_{\mb{x},k}$ is related to $\mb{x}_{k-1}$, so 
\begin{equation}
\begin{aligned}
\mb{A}_{k, k-1, \mb{x}, \mb{x}} &= \frac{\partial \mb{r}_{\mb{x},k}}{\partial \mb{x}_{k-1}} \\
&= - \left( \mathbf{I} - \mathbf{K}_k \mb{H}_k \right) \frac{\partial \mb{x}_{k}}{\partial \mb{x}_{k-1}} \\
&= - \left( \mathbf{I} - \mathbf{K}_k \mb{H}_k \right) \mathbf{F}_k,
\end{aligned}
\end{equation}
according to Eq.~\eqref{eq:KF_residual}, $\mb{r}_{\mb{x},k}$ is related to $\mathbf{K}_{k}$, and $\mathbf{K}_{k}$ is related to $\mathbf{P}_{k-1}$, so
\begin{equation*}
\mb{A}_{k, k-1, \mb{x}, \mathbf{P}} = \frac{\partial \mb{r}_{\mb{x},k}}{\partial \mathbf{P}_{k-1}}
= -\frac{\partial \mathbf{K}_{k}}{\partial \mathbf{P}_{k-1}} \left(\mb{z}_{k} - \mb{H}_{k} \mb{x}_{k-1}\right).
\end{equation*}
Assume that $\mb{S}_k \triangleq \mb{H}_{k}\mathbf{P}_{k-1}\mb{H}_{k}^\intercal + \mb{R}_{k}$. We can further deduce the following results:
\begin{equation*}
\frac{\partial \mathbf{K}_{k}}{\partial \mathbf{P}_{k-1}}
= \frac{\partial \mathbf{K}_{k}}{\partial \mathbf{P}_{k|k-1}}
\frac{\partial \mathbf{P}_{k|k-1}}{\partial \mathbf{P}_{k-1}},
\end{equation*}
first we compute the first term ${\partial \mathbf{K}_{k}}/{\partial \mathbf{P}_{k|k-1}}$,
\begin{equation*}
\begin{aligned}
\d\mathbf{K}_{k}
&= (\d\mathbf{P}_{k|k-1}) \mb{H}_{k}^\intercal \mb{S}_{k}^{-1}
- \mathbf{P}_{k-1} \mb{H}_{k}^\intercal \mb{S}_{k}^{-1} \mb{H}_{k} (\d\mathbf{P}_{k|k-1}) \mb{H}_{k}^\intercal \mb{S}_{k}^{-1}, \\
\mathrm{vec}(\d\mathbf{K}_{k}) 
&= \mathrm{vec}((\d\mathbf{P}_{k|k-1}) \mb{H}_{k}^\intercal \mb{S}_{k}^{-1})
- \mathrm{vec}(\mathbf{P}_{k|k-1} \mb{H}_{k}^\intercal \mb{S}_{k}^{-1} \mb{H}_{k} (\d\mathbf{P}_{k|k-1}) \mb{H}_{k}^\intercal \mb{S}_{k}^{-1}) \\
\mathrm{vec}(\d\mathbf{K}_{k})
&= ((\mb{H}_{k}^\intercal \mb{S}_{k}^{-1})^\intercal \otimes \mathbf{I}) \, \mathrm{vec}(\d\mathbf{P}_{k|k-1}) \\
&\quad - ((\mb{H}_{k}^\intercal \mb{S}_{k}^{-1})^\intercal \otimes \mathbf{P}_{k|k-1} \mb{H}_{k}^\intercal \mb{S}_{k}^{-1} \mb{H}_{k}) \, \mathrm{vec}(\d\mathbf{P}_{k|k-1}) \\
\mathrm{vec}(\d\mathbf{K}_{k}) 
&= \left( \mb{S}_{k}^{-\intercal} \mb{H}_{k} \otimes \mathbf{I}
- \mb{S}_{k}^{-\intercal} \mb{H}_{k} \otimes (\mathbf{P}_{k|k-1} \mb{H}_{k}^\intercal \mb{S}_{k}^{-1} \mb{H}_{k}) \right) \mathrm{vec}(\d\mathbf{P}_{k|k-1}) \\
\mathrm{vec}(\d\mathbf{K}_{k}) 
&= \left( \mb{S}_{k}^{-\intercal} \mb{H}_{k} \otimes (\mathbf{I} - \mathbf{K}_{k} \mb{H}_{k}) \right) \mathrm{vec}(\d\mathbf{P}_{k|k-1}) \\
\frac{\partial \mathbf{K}_{k}}{\partial \mathbf{P}_{k|k-1}}
&= \mb{S}_{k}^{-\intercal} \mb{H}_{k} \otimes (\mathbf{I} - \mathbf{K}_{k} \mb{H}_{k}).
\end{aligned}
\end{equation*}
then we compute the second term ${\partial \mathbf{P}_{k|k-1}}/{\partial \mathbf{P}_{k-1}}$ according to Eq.~\eqref{eq:Kalman_predict}:
\begin{equation*}
\frac{\partial \mathbf{P}_{k|k-1}}{\partial \mathbf{P}_{k-1}} = \mathbf{F}_{k} \otimes \mathbf{F}_{k},
\end{equation*}
thus, the result for $\mb{A}_{k, k-1, \mb{x}, \mathbf{P}}$ is
\begin{equation}
\mb{A}_{k, k-1, \mb{x}, \mathbf{P}}
= -\left[ \left( \mb{S}_{k}^{-\intercal} \mb{H}_{k} \right) \otimes \left( \mathbf{I} - \mathbf{K}_{k} \mb{H}_{k} \right) \right]
\left( \mathbf{F}_{k} \otimes \mathbf{F}_{k} \right) \left( \mb{z}_{k} - \mb{H}_{k} \mb{x}_{k-1} \right),
\end{equation}
according to Eq.~\eqref{eq:KF_residual}, $\mb{R}_{\mb{P},k}$ is not related to $\mb{x}_{k-1}$, so
\begin{equation}
\mb{A}_{k, k-1, \mathbf{P}, \mb{x}} = \frac{\partial \mb{R}_{\mb{P},k}}{\partial \mb{x}_{k-1}} = \mb{0},
\end{equation}
according to Eq.~\eqref{eq:KF_residual}, $\mb{R}_{\mb{P},k}$ is related to $\mathbf{P}_{k-1}$, so
\begin{equation*}
\begin{aligned}
\d\mb{R}_{\mb{P},k} &= (\d\mathbf{K}_{k}) \mb{H}_{k} \mathbf{P}_{k-1}
- (\mathbf{I} - \mathbf{K}_{k} \mb{H}_{k}) (\d\mathbf{P}_{k-1}) \\
\d\mb{R}_{\mb{P},k} &= (\d\mathbf{P}_{k-1}) \mb{H}_{k}^\intercal \mb{S}_{k}^{-1} \mb{H}_{k} \mathbf{P}_{k-1}
- \mathbf{P}_{k-1} \mb{H}_{k}^\intercal \mb{S}_{k}^{-1} \mb{H}_{k} (\d\mathbf{P}_{k-1}) \mb{H}_{k}^\intercal \mb{S}_{k}^{-1} \mb{H}_{k} \mathbf{P}_{k-1} \\
& \quad - (\mathbf{I} - \mathbf{K}_{k} \mb{H}_{k}) (\d\mathbf{P}_{k-1}) \\
\d\mb{R}_{\mb{P},k} &= \mathbf{F}_{k} (\d\mathbf{P}_{k-1}) \mathbf{F}_{k}^\intercal \mb{H}_{k}^\intercal \mb{S}_{k}^{-1} \mb{H}_{k} \mathbf{P}_{k-1} \\
& \quad - \mathbf{P}_{k-1} \mb{H}_{k}^\intercal \mb{S}_{k}^{-1} \mb{H}_{k} \mathbf{F}_{k} (\d\mathbf{P}_{k-1}) \mathbf{F}_{k}^\intercal \mb{H}_{k}^\intercal \mb{S}_{k}^{-1} \mb{H}_{k} \mathbf{P}_{k-1} \\
& \quad - (\mathbf{I} - \mathbf{K}_{k} \mb{H}_{k}) \mathbf{F}_{k} (\d\mathbf{P}_{k-1}) \mathbf{F}_{k}^\intercal,
\end{aligned}
\end{equation*}
\begin{equation*}
\begin{aligned}
\mathrm{vec}(\d\mb{R}_{\mb{P},k}) &= \mathrm{vec}\left(\mathbf{F}_{k} (\d\mathbf{P}_{k-1}) \mathbf{F}_{k}^\intercal \mb{H}_{k}^\intercal \mb{S}_{k}^{-1} \mb{H}_{k} \mathbf{P}_{k-1}\right) \\
& \quad - \mathrm{vec}\left(\mathbf{P}_{k-1} \mb{H}_{k}^\intercal \mb{S}_{k}^{-1} \mb{H}_{k} \mathbf{F}_{k} (\d\mathbf{P}_{k-1}) \mathbf{F}_{k}^\intercal \mb{H}_{k}^\intercal \mb{S}_{k}^{-1} \mb{H}_{k} \mathbf{P}_{k-1}\right) \\
& \quad - \mathrm{vec}\left((\mathbf{I} - \mathbf{K}_{k} \mb{H}_{k}) \mathbf{F}_{k} (\d\mathbf{P}_{k-1}) \mathbf{F}_{k}^\intercal\right) \\
&= \left( (\mathbf{F}_{k}^\intercal \mb{H}_{k}^\intercal \mb{S}_{k}^{-1} \mb{H}_{k} \mathbf{P}_{k-1})^\intercal \otimes \mathbf{F}_{k} \right) \mathrm{vec}(\d\mathbf{P}_{k-1}) \\
& \quad - \left( (\mathbf{F}_{k}^\intercal \mb{H}_{k}^\intercal \mb{S}_{k}^{-1} \mb{H}_{k} \mathbf{P}_{k-1})^\intercal \otimes (\mathbf{P}_{k-1} \mb{H}_{k}^\intercal \mb{S}_{k}^{-1} \mb{H}_{k} \mathbf{F}_{k}) \right) \mathrm{vec}(\d\mathbf{P}_{k-1}) \\
& \quad - (\mathbf{F}_{k} \otimes ((\mathbf{I} - \mathbf{K}_{k} \mb{H}_{k}) \mathbf{F}_{k})) \mathrm{vec}(\d\mathbf{P}_{k-1}) \\
&= \left( (\mathbf{P}_{k-1}^\intercal \mb{H}_{k}^\intercal \mb{S}_{k}^{-\intercal} \mb{H}_{k} \mathbf{F}_{k}) \otimes ((\mathbf{I} - \mathbf{K}_{k} \mb{H}_{k}) \mathbf{F}_{k}) - \mathbf{F}_{k} \otimes ((\mathbf{I} - \mathbf{K}_{k} \mb{H}_{k}) \mathbf{F}_{k}) \right) \mathrm{vec}(\d\mathbf{P}_{k-1}) \\
&= \left( ((\mathbf{P}_{k-1}^\intercal \mb{H}_{k}^\intercal \mb{S}_{k}^{-\intercal} \mb{H}_{k} - \mathbf{I}) \mathbf{F}_{k}) \otimes ((\mathbf{I} - \mathbf{K}_{k} \mb{H}_{k}) \mathbf{F}_{k}) \right) \mathrm{vec}(\d\mathbf{P}_{k-1}),
\end{aligned}
\end{equation*}
thus, the result for $\mb{A}_{k, k-1, \mathbf{P}, \mathbf{P}}$ is
\begin{equation}
\mb{A}_{k, k-1, \mathbf{P}, \mathbf{P}}
= \left[ \left( \mathbf{P}_{k-1}^\intercal \mb{H}_{k}^\intercal \mb{S}_{k}^{-\intercal} \mb{H}_{k} - \mathbf{I} \right) \mathbf{F}_{k} \right] 
\otimes 
\left[ (\mathbf{I} - \mathbf{K}_{k} \mb{H}_{k}) \mathbf{F}_{k} \right],
\end{equation}
according to Eq.~\eqref{eq:KF_residual}, $\mb{r}_{\mb{x},k}$ is related to $\mb{x}_{k}$, so
\begin{equation}
\mb{A}_{k, k, \mb{x}, \mb{x}} = \frac{\partial \mb{r}_{\mb{x},k}}{\partial \mb{x}_{k}} = \mathbf{I},
\end{equation}
according to Eq.~\eqref{eq:KF_residual}, $\mb{r}_{\mb{x},k}$ does not depend on $\mathbf{P}_{k}$. 
When differentiating with respect to a matrix variable, we first vectorize the matrix and define
\[
\frac{\partial\,\mathrm{vec}(f(\mathbf{P}))}{\partial\,\mathrm{vec}(\mathbf{P})},
\]
since $\mb{r}_{\mb{x},k}$ contains no $\mathbf{P}_{k}$, any perturbation $\Delta\mathbf{P}_{k}$ produces zero change, so
\begin{equation}
\mb{A}_{k, k, \mb{x}, \mathbf{P}} 
= \frac{\partial \mb{r}_{\mb{x},k}}{\partial \mathbf{P}_{k}} 
= \mb{0},
\end{equation}
according to Eq.~\eqref{eq:KF_residual}, $\mb{R}_{\mb{P},k}$ is not related to $\mb{x}_{k}$, so
\begin{equation}
\mb{A}_{k, k, \mathbf{P}, \mb{x}} = \frac{\partial \mb{R}_{\mb{P},k}}{\partial \mb{x}_{k}} = \mb{0},
\end{equation}
similarly, $\mb{R}_{\mb{P},k}$ depends linearly on $\mathbf{P}_{k}$ (essentially an identity mapping). 
Thus, for any perturbation $\Delta\mathbf{P}_{k}$ we have 
$\mathrm{d}\mb{R}_{\mb{P},k} = \Delta\mathbf{P}_{k}$, 
which under vectorization gives the identity operator. Hence
\begin{equation}
\mb{A}_{k, k, \mathbf{P}, \mathbf{P}} 
= \frac{\partial \mb{R}_{\mb{P},k}}{\partial \mathbf{P}_{k}} 
= \mathbf{I},
\end{equation}
this confirms that the covariance residual with respect to $\mathbf{P}_{k}$ simply yields the identity operator.

To summarize, the block components of $\mb{A}_{k,k-1}$ and $\mb{A}_{k,k}$ are given by
\begin{align}
\mb{A}_{k, k-1, \mb{x}, \mb{x}} &= - \left( \mathbf{I} - \mathbf{K}_{k} \mb{H}_{k} \right) \mathbf{F}_{k} \label{eq:Akk-1xx} \\
\mb{A}_{k, k-1, \mb{x}, \mathbf{P}} &= - \left[ \mb{S}_{k}^{-\intercal} \mb{H}_{k} \otimes (\mathbf{I} - \mathbf{K}_{k} \mb{H}_{k}) \right] 
(\mathbf{F}_{k} \otimes \mathbf{F}_{k}) (\mb{z}_{k} - \mb{H}_{k} \mb{x}_{k-1}) \label{eq:Akk-1xP} \\
\mb{A}_{k, k-1, \mathbf{P}, \mb{x}} &= \mb{0} \label{eq:Akk-1Px} \\
\mb{A}_{k, k-1, \mathbf{P}, \mathbf{P}} &= \left[ ( \mathbf{P}_{k-1}^\intercal \mb{H}_{k}^\intercal \mb{S}_{k}^{-\intercal} \mb{H}_{k} - \mathbf{I} ) \mathbf{F}_{k} \right]
\otimes \left[ (\mathbf{I} - \mathbf{K}_{k} \mb{H}_{k}) \mathbf{F}_{k} \right] \label{eq:Akk-1PP} \\
\mb{A}_{k, k, \mb{x}, \mb{x}} &= \mathbf{I} \label{eq:Akkxx} \\
\mb{A}_{k, k, \mb{x}, \mathbf{P}} &= \mb{0} \label{eq:AkkxP} \\
\mb{A}_{k, k, \mathbf{P}, \mb{x}} &= \mb{0} \label{eq:AkkPx} \\
\mb{A}_{k, k, \mathbf{P}, \mathbf{P}} &= \mathbf{I} \label{eq:AkkPP}
\end{align}
these expressions complete the characterization of the Jacobian blocks used in our formulation.


\section{Verification of sensitivity}
\label{sec:Verification of sensitivity}

In order to verify the correctness of the above block Jacobian matrix, we first calculate the results of central finite difference and automatic differentiation, and then compare them with the derived theoretical solution. 
In terms of the dynamics and initial parameter setting, we use the same as Eq.~\eqref{eq:rocket_dynamics}.

We set the number of time steps to 100 and plot the Frobenius error norms of finite difference and automatic differentiation, as well as the Frobenius error norms of the theoretical solution and finite difference for the above non-zero Jacobian matrix. In the plots, ``TH'' represents the formula results that we have derived, and ``AD'' and ``FD'' represent automatic differentiation results and finite difference results. From the numerical results, we observe that the error between the theoretical and finite difference (‖TH–FD‖) remains consistently small (typically below $10^{-5}$), and the error between the automatic differentiation and finite difference (‖AD–FD‖) is also at a similar level. This demonstrates that the analytical Jacobian, the automatic differentiation, and the finite difference approximations are all mutually consistent and accurate.

\begin{figure}[htbp]
\centering

\begin{subfigure}[b]{1.\textwidth}
\centering
\includegraphics[width=\textwidth]{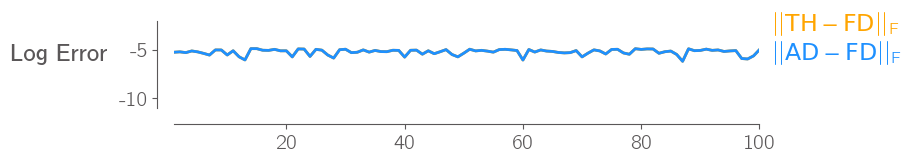}
\label{fig:err_r_x_xk}
\end{subfigure}

\begin{subfigure}[b]{1.\textwidth}
\centering
\includegraphics[width=\textwidth]{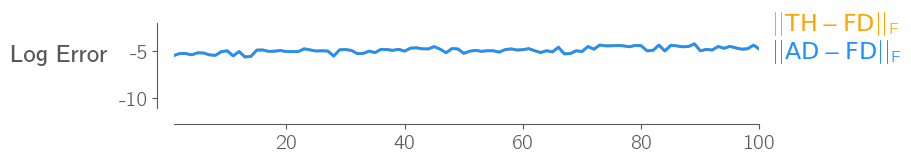}
\label{fig:err_r_x_xprev}
\end{subfigure}

\begin{subfigure}[b]{1.\textwidth}
\centering
\includegraphics[width=\textwidth]{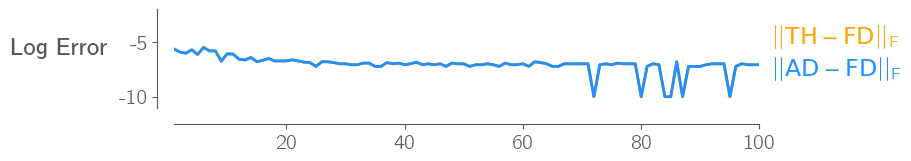}
\label{fig:err_R_P_Pk}
\end{subfigure}

\begin{subfigure}[b]{1.\textwidth}
\centering
\includegraphics[width=\textwidth]{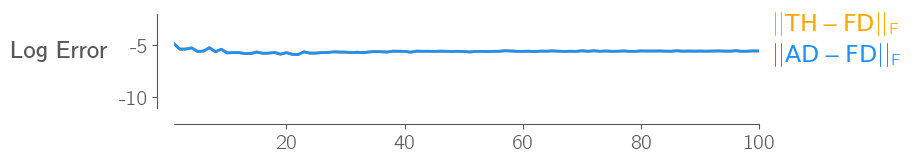}
\label{fig:err_r_x_Pprev}
\end{subfigure}

\begin{subfigure}[b]{1.\textwidth}
\centering
\includegraphics[width=\textwidth]{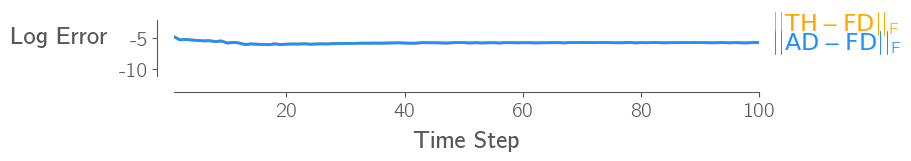}
\label{fig:err_R_P_Pprev}
\end{subfigure}

\caption{F-norm error between automatic differentiation (AD), finite difference (FD), 
and theoretical Jacobian (TH) for various partial derivatives across time steps. 
From top to bottom: 
Eqs.~\eqref{eq:Akkxx}, 
~\eqref{eq:Akk-1xx}, 
~\eqref{eq:AkkPP}, 
~\eqref{eq:Akk-1xP}, 
~\eqref{eq:Akk-1PP}.}
\label{fig:jacobian_error_all}
\end{figure}



\end{document}